\documentclass{amsart}
\usepackage{amssymb}   
\usepackage{amsmath}
\usepackage{amsthm}
\def\a={{\buildrel a \over =}}

\def\GG{{\mathbb G}}

\def\PP{{\mathbb P}}

\def\spr{{\rm Sym}^2}

\newcommand\proofsquare{\nobreak\hfill \hbox{%
\vrule height 5pt
\kern-.4pt
 \vbox{%
\hrule width 5pt depth0pt height.4pt
 \kern4.6pt \hrule  }
\kern-3.75pt
\vrule height 5pt}\kern1pt
\par}

\newtheorem{theorem}{Theorem}[section]
\newtheorem{lemma}[theorem]{Lemma}
\newtheorem{proposition}[theorem]{Proposition}
\newtheorem{corollary}[theorem]{Corollary}

\newtheorem{definition-lemma}[theorem]{Definition-Lemma}

\theoremstyle{definition}
\newtheorem{definition}[theorem]{\bf Definition}

\newtheorem{notation}[theorem]{\bf Notation}

\theoremstyle{remark}
\newtheorem{remark}[theorem]{\bf Remark}

\newtheorem{inters}[theorem]{\bf Inters}

\begin{document}
\title[Note on Fano Surfaces of Nodal Cubic Threefolds]
{A Note on Fano Surfaces of Nodal Cubic Threefolds}
\author{Gerard van der Geer $\,$}
\address{Korteweg-de Vries Instituut, Universiteit van
Amsterdam, Plantage Muidergracht 24, 1018 TV Amsterdam, The Netherlands}
\email{geer@science.uva.nl}
\author{$\,$ Alexis Kouvidakis}
\address{Department of Mathematics, University of Crete, GR-71409 Heraklion, Greece}
\email{kouvid@math.uoc.gr}

\begin{abstract}
We study the Picard variety of the Fano surface of nodal and mildly cuspidal
cubic threefolds in arbitrary characteristic 
by relating divisors on the Fano surface to divisors on
the symmetric product of a curve of genus~$4$.
\end{abstract}
\maketitle
\begin{section}{Introduction}
Cubic threefolds have been studied extensively, first by the classical geometers
starting with Fano \cite{F} and later by Clemens and Griffiths \cite{CG} 
and many others. Nodal cubic threefolds were already considered by Clemens and
Griffiths. In their paper the intermediate Jacobian plays a central role. 
In this note we  come back to these nodal threefolds, but we let 
the Picard variety of the Fano surface of lines on the cubic threefolds 
replace the intermediate Jacobian. 
This has the advantage that it works in all characteristics,
including characteristic $2$. 
We relate divisors on the Fano surface directly to curves on the symmetric 
product of a non-hyperelliptic curve of genus $4$.

We begin by showing that for a canonically embedded non-hyperelliptic
curve $C \subset {\PP}^3$ of genus $4$ the linear system of cubics 
passing through $C$ maps ${\PP}^3$ to a cubic threefold with a node or a 
mild cusp.
The Fano surface of lines on such a cubic threefold is a non-normal 
surface whose normalization equals $\spr C$. We analyze for
linear systems on a curve  $C$ the associated  trace divisors on 
$\spr C$  and their intersection theory
and apply this to divisors on the Fano surface. We analyze
the Picard variety of the Fano surface. 
As an application we get a variation of the proof by Collino 
that the general cubic threefold is not rational and this variation
works also in characteristic~$2$. 
We also analyze the standard compactification of the
Picard variety of the Fano surface, compare the Clemens-Griffiths
map $S \to {\rm Pic}^0(S)$ to the Abel-Jacobi map $\spr C \to {\rm Pic}^0(C)$
and derive a formula for the algebraic equivalence class of 
the Abel-Jacobi image of the Fano surface.

We work over an algebraically closed field $k$ of arbitrary characteristic.
\end{section}
\begin{section}{The Nodal and Mildly Cuspidal Cubic Threefold}\label{nct}
We start with a non-hyperelliptic curve $C$ of genus $4$. The canonical
map $\phi_K: C\to  {\mathbb P}^3$ has as image a curve of degree $6$.
There is a unique quadric $Q$ in ${\mathbb P}^3$ passing through $\phi_K(C)$.
If $Q$ is smooth then the two rulings of $Q$ determine two $1$-dimensional
linear systems $|D_1|$ and $|D_2|$ on $C$ of degree $3$ 
($g^1_3$'s) with $D_1+D_2\sim K$. 
This determines two embeddings $\gamma_i: C \to\spr C $
of $C$ into the symmetric product of $C$, 
by  sending $P \in C$ to the effective divisor of degree $2$ in $|D_i-P|$. 
The images $\gamma_1(C)$ and $\gamma_2(C)$ are disjoint.
If $Q$ is singular then it is of rank $3$ and the two $g^1_3$'s on $C$ coincide
and we find only one embedding $\gamma: C \to\spr C $
sending $P$ to the pair $P_1+P_2$ such that $P+P_1+P_2$
is a divisor of the unique $g^1_3$.

There is a $4$-dimensional linear system $\Pi$ of
cubic surfaces passing through  $\phi_K(C)$.
This defines a rational map
$$
\rho: {\mathbb P}^3\to {\mathbb P}^4=\Pi^{\vee},
$$
sending a point $p\in {\mathbb P}^3-\phi_K(C)$
to the hyperplane $\{ H \in \Pi:  p \in H\}$.
Note that $\Pi$ is the projective space of the kernel 
$U$ of the surjective linear map
${\rm Sym}^3(H^0(C,\omega_C)) \to H^0(C,\omega_C^{\otimes 3})$,
where $\omega_C$ is the relative dualizing sheaf of $C$.
The image is a cubic hypersurface $X$ in ${\mathbb P}^4$.
To see this, one may restrict $\rho$ to a general plane $T$ in 
${\mathbb P}^3$ and obtain the map given by the cubics
passing through six general points, the intersection points
of $\phi_K(C)$ with $T$. The image of $T$ is a Del Pezzo
surface of degree $3$ and this is a linear section of $X$.

We need a precise description of the map $\rho$. Observe that a cubic
hypersurface containing $\phi_K(C)$ that contains also a point
of $P\in Q-\phi_K(C)$ automatically contains $Q$ for reasons
of degree. Therefore the image of such a point
is independent of the choice of $P\in Q-\phi_K(C)$ and 
the open part $Q-\phi_K(C)$ of the quadric $Q$ is contracted. 
But $\phi_K(C)$ is blown up to a ${\mathbb P}^1$-bundle,
the projectivized normal bundle of $\phi_K(C)$. Thus
the image $X$ can be obtained by first blowing up 
${\mathbb P}^3$ along $\phi_K(C)$
and then blowing down the proper transform $Q'$ of $Q$. 
So $X$ has one singular point $x_0$, the image of $Q'$. If $Q$ is smooth then the point $x_0$ is
a node singularity (type $A_1$). If  $Q$ is singular then the point $x_0$ is
a cusp singularity (type $A_2$ if ${\rm char}(k)\neq 2$).
We therefore have the diagram
\[
\begin{matrix}
\tilde{X} && {\buildrel \pi \over\longrightarrow} && X &\subset {\PP}^4 \cr
&\searrow{b} && \nearrow{\rho}   \cr
&& {\mathbb P}^3 \cr
\end{matrix}
\]
where $\tilde{X}$ is the proper transform of ${\mathbb P}^3$
under the blow up map $b$ along $\phi_K(C)$ and
$\pi$ is the blow down map of the proper transform $Q'$ of the quadric surface $Q$.
       
A point $P_1+P_2$ of $\spr C$ determines a line
$\ell \subset {\mathbb P}^3$, namely the line connecting $\phi_K(P_1)$
and $\phi_K(P_2)$ and we interpret this as the tangent line
at $\phi_K(P_i)$ to $\phi_K(C)$ if $P_1=P_2$. 
This line is contained in $Q$ if and only if  
$(P_1,P_2) \in  \gamma_1(C) \cup \gamma_2(C)$ (resp.\
$\gamma(C)$).

Next we need a precise description of the lines in 
${\mathbb P}^4=\Pi^{\vee}$, in particular those that lie on $X$. 
A line in ${\mathbb P}^4$  can be given by a three-dimensional subspace
of the $5$-dimensional space $U$.

A chord $L$ of $\phi_K(C)$ that is not contained in $Q$
determines a line in ${\mathbb P}^4$, that is, the subspace $W$ of $U$
consisting of elements vanishing in two points of the chord 
different from the
intersection points of $L$ with $\phi_K(C)$. Since the elements of $U$
vanish already in $L \cap \phi_K(C)$ this means that the elements of $W$
vanish on $L$. The linear system $\Pi$ restricted to $L$ has projective
dimension $1$, hence maps the chord to a line in ${\mathbb P}^4$.
This line is contained in $X$, and does not pass through $x_0$.

The other lines can be obtained as follows. Take a point $P \in \phi_K(C)$.
The line of a ruling of $Q$ through $P$ determines a point $(P_1,P_2)
\in\spr C $ with $P_1+P_2+P$ a divisor of a $g^1_3$ 
and  a line, namely the image of the exceptional line over $P$ 
in the blow-up of ${\mathbb P}^3$.
Or, differently, the corresponding $3$-dimensional space of $U$ is
the space of cubic surfaces in ${\mathbb P}^3$ 
containing $Q$ and a hyperplane through $P$.
These are lines through the singular point~$x_0$.
So in case $Q$ is smooth the Fano surface $S$ of lines 
is $\spr C$ with $\gamma_1(C)$ and $\gamma_2(C)$ identified.
In case $Q$ is singular the Fano surface has a cusp singularity along 
a curve isomorphic to $C$. Its normalization is $\spr C $ 
and let $\nu: \tilde{S}=\spr C \to  S$ be
the normalization map.

We shall identify  $C$ with its image  $\phi_K(C)\subset {\mathbb P}^3$. 
For $p\in C $ we denote by  $E_p$ the
exceptional line over $p$.
For $s \in S$ we denote by $\ell_s$ the corresponding line in $X$.
To summarize the above description of lines on $X$, we have: 
if $p+q$ is a point of $\spr C$ not on $\gamma_1(C) \cup \gamma_2(C)$
then $\ell_{\nu(p+q)}=\pi_*(\tilde{\overline{pq}})$, the push forward by the map  $\pi $ of 
the proper transform of the secant line $\overline{pq}$. If $p+q \in \gamma_i(C)$
then $\ell_{\nu(p+q)}=\pi_*(E_{\gamma_i^{-1}(p+q)})$.

If $X$ is a cubic threefold with one singularity $x$ which is resolved by
a quadric of rank $3$ such that the tangent cone at $x$ intersects
$X$ along a smooth curve not passing through $x$ then we call $X$ mildly 
cuspidal. 
In \cite{CG} (see also  \cite{CM} and \cite{CML}) 
it is proved that every nodal or mildly cuspidal cubic in ${\PP}^4$
is obtained by blowing up ${\PP}^3$ along a canonically embedded
non-hyperelliptic curve $C$ of genus $4$ and by blowing down the
proper transform of the unique quadric containing $C$. This extends
without problems to characteristic $2$.
We therefore have, cf.\ also \cite{CML} Corollary 3.3:
\begin{proposition}
For a non-hyperelliptic curve $C$ of genus $4$ with canonical image contained 
in a smooth quadric (resp. singular quadric) the linear system of cubics
passing through the canonical image of $C$ maps ${\PP}^3$ to a nodal 
(resp.\ mildly cuspidal)  cubic in ${\PP}^4$. 
This defines an isomorphism between the moduli space 
${\mathcal M}_4-{\mathcal H}_4$ of non-hyperelliptic curves of genus $4$
and the moduli space of nodal or mildly cuspidal cubic threefolds 
in ${\PP}^4$.
\end{proposition} 

\end{section}
\begin{section}{The Symmetric Square of a Curve}\label{symm2}

For a smooth curve  $F$ of genus $g$ and $p\in F$ a point of $F$
we define the divisor $X_p=\{ p+q:\; q \in F \}$  on $ \spr F$
(image of a fiber from the ordinary product) and we denote by 
$x$ its class for algebraic equivalence. If we denote by $j_p: F \to \spr F$ 
the inclusion defined by $j_p(q)=p+q$ then $X_p$ is the isomorphic image of $F$ under $j_p$. 
We shall write
algebraic equivalence as ${\buildrel a \over =}$.
 
We write $\Delta=\{ p+p :\; p \in F\}\subset \spr F$ 
for the diagonal divisor and $\delta $ for its class.  
We have the intersections on $\spr F$: 
$$
x^2=1, \quad x\cdot \delta=2 \quad {\rm and} \quad \delta^2=4(1-g).
$$ 
A divisor (class) $A=\sum_i p_i - \sum_jq_j$ on $F$ 
defines a divisor (class) $S_A=\sum_i X_{p_i} - \sum_jX_{q_j}$ 
on  $\spr F$. The map $A \mapsto S_A$ is obviously linear in $A$. 

\begin{lemma}
\label{Picsym}
For a smooth curve $F$  of genus $g$ the map $O(A) \mapsto O(S_A)$ 
defines an isomorphism
$i: {\rm Pic}^0(F) {\buildrel \sim \over \longrightarrow}   
{\rm Pic}^0(\spr F) $ with inverse map the  $j_p^*: {\rm Pic}^0(\spr F)
\to  {\rm Pic}^0(F)$ for a fixed  $p\in F$.   
\end{lemma}
\begin{proof} We have $j^*_p  \circ  i=1$ 
showing that $i$ is an injection. To show that $j^*_p$ is the inverse map to  $i$ 
it suffices to prove that $i$ is onto. 
The result follows from the fact that 
${\rm Pic}^0(\spr F)$ and ${\rm Pic}^0(F)$ have the same dimension.
This follows from the fact that $\dim H^1(\spr F,O_{\spr F})=g$
as a standard calculation shows.

\end{proof}

\begin{definition}
For a linear system $\Gamma$ of degree $d$ and rank $1$ (a $g^1_d$), 
we define a {\em trace divisor} on $\spr F$ of pairs
contained in $\Gamma$ by  
$$
T_{\Gamma}= \{p+q:\; \Gamma \geq p+q \}.
$$
\end{definition}
We shall denote linear equivalence by $\sim$. Then we have:
\begin{lemma}
\label{basiclemma} 
If $A$ is a degree $d$ divisor on $F$ 
defining a linear system  $\Gamma=|A|$ which is a $g^1_d$ then 
in ${\rm Pic}(\spr F)$ we have 
$$
\Delta \sim 2S_A -2 T_{\Gamma}.
$$
In particular, the divisor class $\Delta$ is divisible by $2$ in
${\rm Pic}(\spr F)$, i.e.\ there is a divisor class $\Delta/2$
such that
$$
T_{\Gamma} \sim  S_A -\Delta/2.
$$
\end{lemma}
\noindent
\begin{proof}
It suffices to prove the result in the case where $\Gamma$ is base 
point free $g^1_d$. Indeed, otherwise the elements of the $g^1_d$ 
have the form $D+D_0$, where $D$ is an element of a base point free 
$g^1_{d-n}$ and $D_0=\sum_{i=1}^n p_i$ is the base divisor. 
Then $A \sim  A' + D_0$, where $A'$ is the divisor of the $g^1_{d-n}$.  
Therefore the points of $T_{\Gamma}$ 
have the form i)  points $a+b$ with $D\geq a+b$ for some  $D\in  \Gamma'$,
the $g^1_{d-n}$, 
and ii) points $a+p_i$, some $a\in F$ and some $i$ with $1\leq i \leq n$. 
The first are points of $T_{\Gamma'}$ while the latter are 
points of $X_{p_i}$,  $i=1,\ldots,n$. We conclude that 
$T_{\Gamma }=T_{\Gamma '} + \sum_{i=1}^nX_{p_i}$ 
and since  $S_A \sim S_{A'}+ \sum_{i=1}^n X_{p_i}$,
the result for $\Gamma'$ implies it for $\Gamma$.

Assume therefore that $\Gamma$ is base point free $g^1_d$
and let $\phi: F \rightarrow {\mathbb P}^1$ be the map defined 
by the $g^1_d$. Take 
$\Phi =\phi \times \phi : 
F\times F \rightarrow {\mathbb P}^1\times {\mathbb P}^1$. 
Let $\beta: F\times F \rightarrow  \spr F$ be the canonical map. 
We denote by $\Delta_{{\mathbb P}^1 \times {\mathbb P}^1}$ 
(resp.\ $\Delta_{F\times F}$) the diagonal of 
${\mathbb P}^1 \times {\mathbb P}^1$ (resp.\ $F\times F$). 
Then  $\Phi ^* \Delta_{{\mathbb P}^1 \times {\mathbb P}^1} 
= \Delta_{F \times F}+\beta^* T_{\Gamma}$.
Now, $\Delta_{{\mathbb P}^1 \times {\mathbb P}^1}$ is linearly
equivalent to $f^1+f^2$, 
where $f^i$, $i=1,2$, are the fibers of ${\PP}^1 \times {\PP}^1$. 
Let $A=a_1+\cdots +a_d$ be a fiber divisor of the map 
$\phi: F \to {\PP}_1$.  
Then  $\Phi ^*\Delta_{{\PP}^1 \times {\PP}^1} \sim  
\Phi ^* (f^1+f^2)= \sum_{i=1}^d\beta^* X_{a_i}=\beta^*S_A$. 
We also have $\beta^*\Delta = 2 \Delta_{F\times F}$. Therefore
$$
\beta^* 2 T_{\Gamma} = 2\Phi ^* \Delta_{{\PP}^1 \times {\PP}^1} 
-2\Delta_{F \times F} 
\sim 2 \beta ^*  S_A- \beta^* \Delta.
$$
Since $\beta^*$ is an injection the result follows.
\end{proof}

\begin{remark}
\label{delta}
{\rm Among all divisor classes $A$ on $\spr F$ with 
$A \, {\buildrel a \over =}\, \delta /2$, the above defined 
divisor class $\Delta/2$ has the characteristic property that
$j_p^*(\Delta/2)=p$ for every $p \in F$. Indeed, one has
$j^*_p(\Delta/2)=j^*_p(S_A-T_{\Gamma})=A-(A-p)=p$.  }
\end{remark}
\end{section}
\begin{section}{Divisors on $\spr C$}
\label{divisors on symC}
 We shall
assume for the rest of the paper  that $Q$ is smooth, or equivalently, that the curve $C$ has
two different $g^1_3$'s, say $|D_1|$ and $|D_2|$ with $D_1+D_2\sim K$, the
canonical divisor of $C$. The corresponding cubic threefold $X$ is then a nodal cubic threefold, 
see Section \ref{nct}. 
Let $R_1$ and $R_2$ be the two rulings of $Q$. If $p\in C$ is a point
and $\ell_i$ a line in the ruling $R_i$ through $p$ then $\ell_i$ 
cuts out on $C$ a divisor $p+p_i+q_i$ in one of the two $g^1_3$'s.
The map $\gamma_i$,  see Section \ref{nct},  sends $p$ to $p_i+q_i \in\spr C $. 
We shall write $C_i$ for the image curve
$\gamma_i(C)$ on $\spr C$.
The map $\gamma_2\gamma_1^{-1}: C_1\to C_2$ sends
$p_1+q_1 \in\spr C $ to the  complementary point $p_2+q_2$.
 For
complementary points $p_1+q_1\in C_1$ and $p_2+q_2\in C_2$ we
have $\ell_{\nu(p_1+q_1)}=\ell_{\nu(p_2+q_2)}$.
Therefore the  normalization map $\nu : \tilde{S}=\spr C \to S$ glues 
the complementary points of the curves $C_1$ and $C_2$. 

We observe now that  $C_i=T_{|D_i|}$, i.e.\ $C_i$ is a trace divisor on $\spr C $, and hence by Lemma \ref{basiclemma} we have 
\begin{equation} 
\label{eqC_i}
C_i \sim S_{D_i} -\Delta/2.
\end{equation}

This relation yields the following corollary, cf.\ also
\cite{Coll}.

\begin{corollary}
The divisors $C_1$ and $C_2$ are algebraically equivalent on
$\spr C$, but not linearly equivalent.
\end{corollary}
\begin{proof} 
We have $S_{D_i}\stackrel{a}{=} 3x$, and hence 
$C_i\stackrel{a}{=}3x-\delta/2$ which proves that $C_1\stackrel{a}{=} C_2$. If we assume that $C_i\sim C_2$ then $S_{D_1}-\Delta/2 \sim S_{D_2}-\Delta/2$, and so $S_{D_1}\sim S_{D_2}$. But then Lemma \ref{Picsym} implies that $D_1\sim D_2$, a contradiction since we assumed that $D_1$ and $D_2$ define  two different $g^1_3$'s.  
\end{proof}
Since the curve $C$ is 
not hyperelliptic we have $h^0(C,O(K-p-q))=2$, for every $p,q, \in C$. We introduce now the following notation:
\begin{notation} For $p+q \in\spr C $ we set 
$D_{p+q}=T_{\Gamma}$ with $\Gamma$ the $g^1_4$ defined by $K-p-q$.
\end{notation}

So $D_{p+q}$ is the trace  divisor on $\spr C$ 
corresponding to the projection of the curve $C$ with center the 
secant line $\overline{pq}$. 
By Lemma \ref{basiclemma}, we have
\begin{equation} \label{O(D(p+q))}
D_{p+q}\sim S_{K-p-q}-\Delta/2 \, .
\end{equation}
If $p+q \notin C_1 \cup C_2$ then  $|K-p-q|$ is base point free 
and defines a $4:1$ map from $C$ to ${\PP}^1$. If $p+q \in C_i$ 
then the linear system $|K-p-q| $ has  the base  point $\gamma_i^{-1}(p+q)$.   

\begin{lemma}
\label{Dpq, p+q in Ci}
If $p+q \in C_1$ then 
$D_{p+q}=C_2+X_{\gamma _1^{-1}(p+q)}$ and, similarly, if
$p+q\in C_2$ then $D_{p+q}=C_1+X_{\gamma _2^{-1}(p+q)}$.
\end{lemma}
\begin{proof}
If $p+q \in C_1$  
then the linear system $|K-p-q| $ has the base point 
$\gamma_1^{-1}(p+q)$.  The elements in $|K-p-q|$ have the form 
$D+ \gamma_1^{-1}(p+q)$, where 
$D\in |K-p-q-\gamma_1^{-1}(p+q)|=|K-D_1|=|D_2|$. 
Then, as in the first paragraph of the proof of Lemma 
\ref{basiclemma}, we have that 
$D_{p+q}=T_{|D_2|}+ X_{\gamma _1^{-1}(p+q)}=C_2+ X_{\gamma _1^{-1}(p+q)}$.
\end{proof}

We now compute several intersection numbers on $\spr C$.

\begin{inters}
\label{Ci dot Cj}
$ [C_i]\cdot [C_j]=(3x-\delta/2)^2=9-6-3=0.$
\end{inters}

\begin{inters}
\label{Xp dot Ci}
$[X_p] \cdot [C_i]=x\cdot(3x-\delta/2)=2$.

The two points of intersection are $p+a$, $p+b$, where $a$ and $b$ are defined by $\gamma _i(p)=a+b$. 
\end{inters}

\begin{inters}
\label{Dpq dot Ci}
 $[D_{p+q}] \cdot [C_i]=(4x-\delta/2) \cdot (3x-\delta/2)=12-7-3=2$.

If $p+q \notin C_1\cup C_2$, the two points of intersection are the 
$\gamma_i(p)$ and $\gamma_i(q)$. Note therefore that, 
in this case,  the divisor $D_{p+q}$ intersects the curves  $C_1$ and $C_2$ 
in complementary points: 
$\gamma_1(p) =\gamma_1 \gamma _2^{-1} (\gamma_2(p))$ and $\gamma_1(q)=\gamma_1 \gamma _2^{-1} (\gamma_2(q))$. This
indicates that the divisor $D_{a+b}, \; a+b \notin C_1\cup C_2$, is the pull back of a Cartier divisor from $S$ - 
we will see this later in a more rigorous way.   
If $p+q\in C_1$ then, by Lemma \ref{Dpq, p+q in Ci}, we have that
$D_{p+q}=C_2+X_{\gamma _1^{-1}(p+q)}$  and the 
points of intersection are the two points of intersection of 
$X_{\gamma _1^{-1}(p+q)}$ with $C_i$, 
that is, the points $p+\gamma^{-1}_1(p+q)$ and $q+\gamma^{-1}_1(p+q)$, 
see Inters \ref{Xp dot Ci}. 
\end{inters}

\begin{inters}
\label{Dpq dot x}
$ D_{p+q} \cdot X_a=(4x-\delta/2)\cdot x= 4-1 = 3.$ 
 
If $p+q \notin C_1\cup C_2$ it corresponds to the three points $a+b_i$, where $b_i$ are the 
three additional  points of intersection with $C$ of the plane defined by the points $p,q$ and $a$. 
If $p+q \in C_1$ then  $D_{p+q}= C_2+ X_{\gamma _1^{-1}(p+q)}$, see Lemma \ref{Dpq, p+q in Ci}, 
and the intersection corresponds to the two points of 
intersection of $C_2$ with $X_a$, see Inters \ref{Xp dot Ci}, plus the point of intersection of 
$X_{\gamma _1^{-1}(p+q)}$ with
$X_a$, i.e.\  the point $ \gamma _1^{-1}(p+q)+a$.  
\end{inters}

\begin{inters}
\label{Dpq dot dp'q'}
$[D_{p+q}]\cdot [D_{p'+q'}]=(4x-\delta/2)^2=16-8-3=5$.

If $p+q, \; p'+q' \notin C_1\cup C_2$ it corresponds to the five 
common secant lines  to the  lines $l_{\nu (p+q)}$  and 
$l_{\nu (p'+q')}$ in $X$, cf. \cite{CG}. If $p+q \in C_1$ then  $D_{p+q}=
C_2+ X_{\gamma _1^{-1}(p+q)}$, see Lemma \ref{Dpq, p+q in Ci}, 
and the intersection corresponds to the sum of the Inters
\ref{Dpq dot Ci} and  Inters \ref{Dpq dot x}.
\end{inters}
\end{section}
\begin{section}{Divisors on the Fano Surface}
\label{divisorFano}

Let $X$ be again a nodal threefold with $S$ the Fano surface.
For each  $s \in S$ we have the divisor 
$$D_s=\{ s'\in S, \; l_{s'} \cap l_s\neq \emptyset \}
$$ 
on $S$ as defined in \cite{CG}. Let $s\in S$ so that $s=\nu (p+q)$ for some $p+q \in\spr C $, where $\nu:
\spr C\to S$ is the normalization map. The following proposition relates 
the divisor $D_s$ on $S$ with the trace divisor $D_{p+q}$  
on $\spr C$.  

\begin{proposition}
\label{nu*Dpq} 
Let ${\rm sing}(S)$ be  the singular locus of $S$ viewed as a Weil divisor
on~$S$. We have 
\begin{enumerate} 
\item{} If $p+q \notin C_1 \cup C_2$ then $D_{\nu(p+q)}= \nu_*D_{p+q} $.
\item{} If $p+q\in C_1$   then $D_{\nu(p+q)}=\nu_*D_{p+q} 
={\rm sing}(S)+\nu_*X_{\gamma_1^{-1}(p+q)}$. Similarly, if $p+q\in C_2$   then 
$D_{\nu(p+q)}= \nu_*D_{p+q}  ={\rm sing}(S)+\nu_*X_{\gamma_2^{-1}(p+q)}$.
\end{enumerate}
\end{proposition}
\begin{proof} 
We start by proving the first claim. 
The points $a,b\in C$ belong to the same fiber of the projection to ${\mathbb P}^1$ defined by the $g^1_4=|K-p-q|$  if and only if there is a hyperplane section $H$ 
on $C$ with $H \geq p+q+a+b$. This is equivalent to saying that the  
line $\overline{ab}$ intersects the line  $\overline{pq}$. If 
$p+q\notin C_1\cup C_2$ then the secant $\overline{pq}$ corresponds, 
via the rational map $\rho$, to the line  $l_{\nu(p+q)}=
 \pi_*\tilde{\overline{pq}}$ of $X$.   
The point $p$ (resp.\ $q$) of $\overline{pq}$ corresponds  
to the  intersection $x_p$ (resp.\ $x_q$) of $l_{\nu(p+q)}$ 
with $\pi_*E_{p}$ (resp. $\pi_*E_{q}$). 

Apart from the line $\pi_*E_p$ ($=l_{\nu(\gamma _1(p))}= l_{\nu(\gamma _2(p))}$), 
the lines in $X$ which intersect the line $l_{\nu(p+q)}$ at $x_p$ are the 
lines  $l_{\nu(p+a)}$,  where $p+a$ is one of  three points of intersection of $D_{p+q}$ with $X_p$, 
see Inters \ref{Dpq dot x}. This is because the proper transform 
of two lines through $p$ passes from the same point of $E_p$ 
if and only if the plane they span contains  the tangent line $T_pC$.  
Let 
$S_p=[X_p \cap D_{p+q}]\cup \{ \gamma _1(p), \gamma _2(p) \} \subset D_{p+q}$. 
Then $a+b \in S_p$ if and only if  $l_{\nu(a+b)}$ is a line in $X$ 
that intersects $l_{\nu(p+q)}$ at $x_p$. 
Note that $p+q \in D_{p+q}$ if and only if the tangents to $C$ 
at $p$ and $q$ are coplanar and in this case $l_{\nu(p+q)}$ 
intersects itself. This gives a characterization of the  
lines of second type, see \cite{CG} Lemma 10.7.
Similarly, the set of $s\in S$ such that the line $l_s$  
intersects the line $l_{\nu(p+q)}$ at  $x_q$ is the image of the set 
$S_q=[X_q \cap D_{p+q}]\cup \{ \gamma _1(q), \gamma _2(q) \} \subset D_{p+q}$.

We set $U=D_{p+q}\backslash [S_p\cup S_q]$ and let $U'$  be the set of points $s$ in the divisor 
$D_{\nu(p+q)}$ such that $l_s$  intersects the line $l_{\nu(p+q)}$ at a 
point different from $x_p$ and $x_q$. 
We shall show that $a+b \in U$ if and only if  $\nu(a+b) \in U'$, which yields the first claim.
We claim that $a+b \in U$ if and only if  the line $\overline{ab}$ 
intersects $\overline{pq}$ at a point $t$  different than 
$p$ and $q$. Indeed, the line $\overline{ab}$ 
intersects $\overline{pq}$ since $a+b \in D_{p+q}$. 
As $a+b \notin X_p+X_q$  we have $\{a,b\} \cap \{p,q\}=\emptyset$ 
and so if we assume that the point of intersection is $p$ or $q$ 
then the line $\overline{ab}$ intersects the curve $C$ at 3 points and hence
it is a line in a ruling. But then, Inters \ref{Dpq dot Ci} yields that 
$a+b \in \{ \gamma_i(p), \gamma_i(q), \;i=1,2 \}$, a contradiction since $a+b \in U$.   
 Hence the lines $l_{\nu(a+b)}$ and $l_{\nu(p+q)}$ are 
intersecting lines with point of intersection the 
$\alpha (t) \neq x_p, x_q$. Therefore $\nu(a+b) \in U'$ and vice versa. 

The second claim follows easily from Lemma \ref{Dpq, p+q in Ci}. 
The curve ${\rm sing}(S)$ corresponds to lines intersecting 
$l_{\nu(p+q)}$ at the singular point of the threefold 
and $\nu_*X_{\gamma_i^{-1}(p+q)}$ corresponds to lines 
intersecting  $l_{\nu(p+q)}$ at the other points.
\end{proof}

Since for every $p+q\in\spr C $ the divisors $D_{p+q}$ have 
algebraic equivalence 
class $4x-\delta/2$, see relation (\ref{O(D(p+q))}),  we have
the following corollary.

\begin{corollary}
\label{alg. equiv}
For every $s\in S$ the divisors $D_s$ are algebraically equivalent. 
\end{corollary}

\begin{remark}
\label{involution}
Note that for $p+q \notin C_1\cup C_2$ the divisor $D_{p+q}$ has an 
involution which sends the point $a+b$ to the residual point in the linear 
system $|K-p-q|$. The induced involution on $D_{\nu(p+q)}$ is the 
one defined in \cite{CG}. 
\end{remark}

For a $2$-plane $V$ in ${\PP}^4$ the set of lines in ${\PP}^4$ meeting
$V$ defines a Cartier divisor $C_V$ on $S$. The corresponding divisor class is the
pull back to $S$ via the natural embedding
$S \to {\rm Gr}(2,5)$ of the natural ample line bundle 
on the Grassmannian.
Let $p_1+q_1 \in\spr C $, but $\notin C_1\cup C_2$ 
and choose a generic plane $H$ in ${\PP}^3$ containing the secant 
$\overline{p_1q_1}$ and intersecting the  curve $C\subset {\PP}^3$ 
in four additional distinct points $p_2,q_2,p_3,q_3$ different from $p_1,q_1$.  
We may assume that  $p_2+q_2, \, p_3+q_3\notin C_1 \cup C_2$ 
and that the lines $\overline{p_iq_i}$, $i=1,2,3$,  
meet at three distinct points not on $C$.  
Therefore, their image under the rational map  $\rho$ are three 
intersecting lines in ${\PP}^4$ which define a 2-plane~$V_0$. 
Note that the rational map  $\rho$ embeds the plane $H$ in a 
hyperplane of ${\PP}^4$ but does not  send it to a 2-plane 
in ${\PP}^4$.   Then the plane $V_0$ intersects $X$ in the sum 
of the three lines  $ \sum_{i=1,2,3}l_{\nu(p_i+q_i)}$ and hence $C_{V_0}=
\sum_{i=1,2,3}D_{\nu(p_i+q_i)}$ is a Cartier divisor on $S$. 
Now, by Proposition \ref{nu*Dpq} and Inters  
\ref{Dpq dot Ci} the divisor 
$D_{\nu(p_1+q_1)}$ intersects the singular locus of $S$ at the divisor  
$A=\nu(\gamma_1(p_1)+\gamma_1(q_1))$, while the divisor 
$D_{\nu(p_2+q_2)}+D_{\nu(p_3+q_3)}$
intersects the singular locus of $S$ at the divisor
$B=\nu(\gamma_1(p_2)+\gamma_1(p_3)+\gamma_1(q_2)+\gamma_1(q_3))$.  Since 
${\rm supp}A \cap {\rm supp}B=\emptyset$,  
the divisor $D_{\nu(p_1+q_1)}$ is a Cartier divisor. 
Hence, if $p+q \notin C_1\cup C_2$ then  the 
divisor $D_{\nu(p+q)}$ defines a line bundle  
${\mathcal O}(D_{\nu(p+q)})$ on the singular surface $S$. 
Combining this with Proposition \ref{nu*Dpq} we have: 

\begin{corollary}
\label{pullback lb}
If $s\in S-{\rm sing}(S)$ so that $s=\nu(p+q)$ with $p+q\notin C_1\cup C_2$,  
then $D_s$ is a Cartier divisor on $S$ and  
${\mathcal O}(D_{p+q})=\nu^* {\mathcal O}(D_s)$.
\end{corollary}

\begin{remark}
\label{not cartier}
If $s \in {\rm sing}(S)$ with  $s=\nu(p+q)$ but $p+q\in C_1\cup C_2$, 
then $D_s={\rm sing}S+\nu_*X_{\gamma_i^{-1}(p+q)}$, $i=1$ or $2$ (see Proposition \ref{nu*Dpq}),
is not a Cartier  divisor on $S$. For example, $\nu_*X_{\gamma_i^{-1}(p+q)}$ 
is not a Cartier divisor since for  $s\in C$ the divisor $X_s$ does not intersects the curves
$C_i$, $i=1,2 $, at complementary points.
\end{remark}
\end{section}
\begin{section}{The Picard Variety of the Fano Surface}\label{PicFano}
We now analyze the Picard variety of the Fano surface of our
nodal cubic threefold.
\begin{proposition}
\label{onto}
The pull back map 
$\nu ^*: {\rm Pic}^0 (S) \to {\rm Pic}^0(\tilde{S})$ is onto. 
\end{proposition} 
\begin{proof}
By Lemma  \ref{Picsym} the group ${\rm Pic}^0(\tilde{S})$ is generated 
by the classes of divisors of the form 
$S_{a-b}$ with $a, b  \in C$. 
Choosing a point $c\in C$ with  $a+c,\,  b+c \notin C_1\cup C_2$ we get
by relation (\ref{O(D(p+q))}) and  Corollary  \ref{pullback lb} that
$S_{a-b}$ is linearly equivalent to
$$
[S_{K-b-c}-\Delta/2]-[S_{K-a-c}-\Delta/2]=
D_{b+c}-D_{a+c}=\nu^*(D_{\nu(b+c)}-D_{\nu(a+c)}).
$$
\end{proof}

\begin{remark}
\label{remonto} A line bundle $L$ on $\spr C$ defining an
element of  ${\rm Pic}^0(\spr C)$ restricts to the same
line bundle on $C_1$ and $C_2$, that is, 
$\gamma_1^*(L)\cong \gamma_2^*(L)$.
Indeed, for $p\in C$ the intersection of $X_p$ with $C_i$ is
$s_i+p$, $t_i+p$, where $\gamma_i(p)=s_i+t_i$,  see Inters \ref{Xp dot Ci}.
Therefore $\gamma_i^*(X_p)=s_i+t_i$ and $\gamma_i^*(S_p)\sim D_i-p$
because $s_i+t_i+p \sim D_i$. So $\gamma_i^*(S_{p-q}) \sim q-p$
and since these divisors generate ${\rm Pic}^0(\spr C)$ the
result follows. Given now $L$ we can glue $L|C_1$ with $L|C_2$ to
obtain a line bundle on $S$ which under $\nu$ pulls back to $L$.
This proves the surjectivity of Proposition \ref{onto} in a different way.
\end{remark}

\begin{corollary} 
\label{exactseqPic}
The semi-abelian variety ${\rm Pic}^0(S)$ is isomorphic 
to the ${\GG}_m$-extension of ${\rm Pic}^0(C)$ given by 
$D_1-D_2$, the difference of the two $g^1_3$'s.
\end{corollary}
\begin{proof}
The kernel of the surjective map $\nu^*: {\rm Pic}^0(S) \to
{\rm Pic}^0(\tilde{S})\stackrel{i}{\cong} {\rm Pic}^0(C)$ is the algebraic torus ${\GG}_m$. 
More precisely, the fibre over $[L] \in {\rm Pic}^0(\tilde{S})$
consists of the isomorphisms $\gamma_1^*(L)\cong \gamma_2^*(L)$.
Note that $\gamma_1^*(L)$ and $\gamma_2^*(L)$ are isomorphic line
bundles on $C$ as the map
$${\rm Pic}^0(C) {\buildrel i \over \to} {\rm Pic}^0(\tilde{S}) 
{\buildrel \gamma_j^* \over \longrightarrow} {\rm Pic}^0(C)
$$
is given on divisors as $\sum n_i p_i \mapsto \sum n_i (D_j-p_i)$
(cf.\  Remark \ref{remonto}), hence by $-1$. Let ${\mathcal L}$
be a universal line bundle on ${\rm Pic}^0(C) \times \tilde{S}$ constructed via the 
Abel-Jacobi map $u: \tilde{S} \to {\rm Pic}^0(C)$ with $u(p+q)=O(p_0+q_0-p-q)$ for fixed $p_0,q_0\in C$. 
It has the properties ${\mathcal L}| [L] \times \tilde{S} = O(S_L)$ and 
${\mathcal L}|{\rm Pic}^0 C\times \{ p+q \} = O(p_0+q_0-p-q)$.
The line bundle $(1\times \gamma_2)^*({\mathcal L})
\otimes (1\times  \gamma_1)^*({\mathcal L})^{-1}$ on ${\rm Pic}^0(C)\times C$ 
is trivial on each fibre $C$ and hence the pull back of a line bundle
on ${\rm Pic}^0(C)$. To determine which one, we can restrict
to a fibre ${\rm Pic}^0(C)\times \{p \}$ and then it is
seen to equal $O(D_1-D_2)$, since $(1 \times \gamma_j)^*({\mathcal L})|{\rm Pic}^0(C)\times \{p \}=
O(p_0+q_0-\gamma_j(p))=O(p_0+q_0-D_j+p)$). Hence the ${\GG}_m$-extension is
obtained by deleting the zero-section from the line bundle
$O(D_1-D_2)$.
\end{proof}

\begin{remark}
Note that we do not require isomorphisms of ${\GG}_m$-extensions to
be the identity on ${\GG}_m$, hence $O(D_1-D_2)$ and $O(D_2-D_1)$
define isomorphic extensions.
\end{remark}

As a corollary we now can deduce that the general cubic threefold is
not rational, cf.\ \cite{Coll2}, but with no assumptions on the characteristic.

\begin{corollary}
The general cubic threefold is not rational.
\end{corollary}
\begin{proof}
Let ${\mathcal X} \to B$ be a cubic threefold over the spectrum of a 
discrete valuation ring such that the generic fibre $X_{\eta}$ is
a smooth cubic threefold and the special fibre $X_s$ is a nodal cubic threefold.
Then the Picard variety of the Fano surface
$S$ of ${\mathcal X}$ is a semi-stable abelian variety ${\mathcal A}$ of
dimension $5$ with generic fibre $A_{\eta}$ a principally polarized abelian variety
and as special fibre $A_s$ a ${\GG}_m$-extension of ${\rm Jac}(C)$ given
by $\pm (D_1-D_2)$. If ${\mathcal X}$ were rational then $A_{\eta}$ 
would be the Jacobian of a curve of compact type. But if $D_1\neq D_2$
then $D_1-D_2$ is not of the form $p-q$ for points $p, q \in C$.
Hence $A_s$ is not a limit of a Jacobian, and thus $A_{\eta}$ 
cannot be a Jacobian.
\end{proof} 
\begin{remark}
If $X$ is mildly cuspidal then ${\rm Pic}^0(S)$ is an extension of 
${\rm Pic}^0(C)$ by an additive group, hence not the Jacobian of a
curve of compact type.
\end{remark}
\end{section}
\begin{section}{The Compactified Jacobian of a One-nodal Curve}
\label{compjac}

Since the Picard group ${\rm Pic}^0(S)$ is very
similar to the Jacobian of a one-nodal curve 
we first review compactifications of the latter. 
References are \cite{M2, M, AK2, AK3, EGK, C}
for example. 

Let $F$ be a curve with one node $x_0$ and let 
$\mu :\tilde{F} \to F$ be the normalization, where
$\tilde{F}$ is a smooth curve of genus $g$ with $\mu^{-1}(x_0)=\{x_1,x_2\}$. 
Then there is a short exact sequence of algebraic groups
$$
1\to {\GG}_m \to {\rm Pic}^0(F) \to {\rm Pic}^0(\tilde{F})\to 0
$$
and the extension class of this semi-abelian variety is given by
the class of $\pm (x_1-x_2) \in {\rm Pic}^0(\tilde{F})/{\pm 1}$, where we
identify ${\rm Pic}^0(\tilde{F})$ with its dual abelian variety.
We shall write
$$
J:={\rm Pic}^0(\tilde{F}), \quad G:={\rm Pic}^0(F).
$$

There are two ways two compactify $G$, one by explicitly constructing
a geometric compactification (the rank-1-compactification), 
the other one by the moduli interpretation.
In the case at hand they lead to the same result. We begin by defining
the compactified Jacobian $G^c$ of $F$ as the moduli space of rank~$1$ 
torsion-free sheaves on $F$ of degree $0$; here the degree is defined by
$\deg(M)=\chi(M)-\chi(O_F)$. It contains $G={\rm Pic}^0(F)$
as an open part. 

The direct construction of the compactified Jacobian as a variety
is obtained as follows. Take the ${\PP}^1$-bundle $P={\PP}(L \oplus O)$
with the projection $q: P\to J$ over $J$ where $L=O(x_1-x_2)$. 

Recall that in order to lift a morphism $\alpha: X \to J$  for a
variety $X$ to a morphism $\tilde{\alpha}: X \to P$ one must give
an invertible sheaf $M$ on $X$ and a surjective map of sheaves 
$\tau: \alpha^*(L\oplus O) \to M$, see \cite{H}, Ch.\ II, Prop.\ 7.12.
The ${\PP}^1$-bundle $P$  contains two effective divisors $T_1$ and $T_2$ 
given by ${\PP}(L\oplus 0)$ and ${\PP}(0\oplus O)$. 
There exist two sections $t_i: {\rm Pic}^0(C)\to P$ ($i=1,2$) of $q$ 
with image $T_i$ with  $t_1$ corresponding to the projection 
$L\oplus O \to L$ and $t_2$ to $L\oplus O \to O$.
By deleting $T_1$ and $T_2$ from $P$ one gets $G$ back.

Since $O(T_i)\otimes O(1)^{-1}$ is trivial on the fibers of $q$
it is the pull back of a line bundle $\Lambda_i$ on $J$.
One determines $\Lambda_1$ by taking the pull back of the
relation $O(T_1)\otimes O(1)^{-1}=q^*\Lambda_1$ under $t_2$ and one finds
$\Lambda_1=O$ and similarly one gets $\Lambda_2=L^{-1}$.
In particular we get
\begin{equation}\label{T1T2}
O(T_1)\cong O(1), \quad O(T_2)\cong O(1)\otimes q^*(L^{-1}).
\end{equation}

The compactification $G^c$ is the non-normal variety obtained
by glueing $T_1$ to $T_2$ over a translation by $x_1-x_2$.
The smooth part can be identified with $G$ and the singular
locus $\Sigma$ with ${\rm Pic}^{-1}(\tilde{F})$ by associating to a 
line bundle $N$ on $\tilde{F}$ of degree~$-1$ the torsion-free 
sheaf $\mu_*(N)$ of rank $1$ on $F$:
$$
\mu_* : {\rm Pic}^{-1}(\tilde{F}) {\buildrel \sim \over
\longrightarrow} \Sigma, \qquad N \mapsto \mu_*(N).
$$
Note that $\chi(\tilde{F})=\chi(F)+1$.

There is a natural action of $J$ by translation on itself and
on ${\rm Pic}^{-1}(\tilde{F})$ and this results in an action of $G$ on $G^c$
extending the action on itself as one easily checks. 

Moreover, for a line bundle $\Lambda$ on $\tilde{F}$
the fibre under $\mu^*$ can be interpreted as the pairs $(\Lambda,\lambda)$
where $\lambda: \Lambda_{x_1} {\buildrel \sim \over \to} \Lambda_{x_2}$.
If we choose generators for the fibres $\Lambda_{x_i}$ the map
$\lambda$ can be identified with a non-zero scalar. Letting
this scalar go to $0$ or infinity gives the two extra points on 
the fibre of $G^c$ over $\Lambda$; these have as their images
in $\Sigma$ the points corresponding to the torsion-free sheaves 
$\mu_*(\Lambda \otimes O(-x_1))$ and $\mu_*(\Lambda\otimes O(-x_2)$.

After choosing a smooth point $p_0$ on $F$ with inverse image $\tilde{p}_0$
on $\tilde{F}$ we can define an Abel-Jacobi map
$$
u: \tilde{F} \to J, \quad  \tilde{p} \to O(\tilde{p}_0-\tilde{p})
$$
and it can be lifted to a map $\tilde{u}:
\tilde{F} \to P$ which is given by an invertible sheaf $M$ on $\tilde{F}$
and a surjection $u^*(L\oplus O) \to M$ with $M=\tilde{u}^*(O(1))$.
We take $M=O(x_1)\cong O(x_2)\otimes u^*(L)$. Then 
$M \otimes u^*(L\oplus O )^{\vee}
\cong O(x_1)\oplus O(x_2)$ and this has a canonical section $1\oplus 1$,
giving the desired surjection $u^*(L\oplus O) \to M$.

Note that by equation 
(\ref{T1T2})  for $i=1,2$ we have $\tilde{u}^*O(T_i)=O(x_i)$ 
as $\tilde{u}^*O(1)=O(x_1)$ and $u^*L^{-1}=O(x_2-x_1)$.
It follows that $\tilde{u}(\tilde{F})=\tilde{u}(x_i)$ and 
$\tilde{F}-\{x_1,x_2\}$ is mapped into $G$ under $\tilde{u}$. 
Note also that $q(\tilde{u}(x_1))-q(\tilde{u}(x_2))=u(x_1)-u(x_2)$, 
the class of $x_1-x_2$. Thus the morphism $\tilde{u}$ descends 
to an Abel-Jacobi map $\bar{u}: F \to G^c$. It
has a moduli interpretation via the direct construction as follows.
The ideal sheaf $I_{\Delta}$ of the diagonal on $F\times F$ is a 
torsion-free sheaf of degree $-1$ for the curve $F\times F 
{\buildrel {\rm pr} \over \longrightarrow} F$, with ${\rm pr}$ the first
projection. Then the sheaf $I_{\Delta} \otimes {\rm pr}^*O(p_0)$
defines the morphism $\bar{u}$; we refer to \cite{EGK}.

We now calculate the algebraic equivalence class of the curve 
$\tilde{u}(\tilde{F})$ in $P$.

\begin{proposition}
The algebraic equivalence class $\gamma$ of the curve $\tilde{u}(\tilde{F})$
in $P$ is given by 
$$
\tilde{u}(\tilde{F}) \, {\buildrel a \over =} \, q^*(p)+q^*(F) \cdot \eta,
$$
with $q: P \to J$ the projection and $\eta$ the class of $O(1)$ on $P$
and $p$ a point of $J$.
\end{proposition}
\begin{proof}
From equation (\ref{T1T2}) we have $T_1\a= \eta$ and $T_2 \a= \eta - q^*(L)\a= \eta$.
Since $\eta^2=0$ we have $\gamma \a= q^*(a_0)+q^*(a_1)\eta$ 
with $a_i$ a class of dimension $i$ on $J$ satisfying
$a_0=q_*(\gamma \eta)$ and $a_1=q_*(\gamma)$. We have
$$
q_*(\gamma)=q_*\tilde{u}_* 1_{\tilde{F}} =u_* 1_F=[\tilde{F}]
$$
and $q_*(\gamma \eta)$ equals the class of a point $p$
and the formula follows.
\end{proof}

\end{section}
\begin{section}{The Compactified Picard of the Fano Surface}
\label{compfano}
As we saw in section \ref{PicFano} the semi-abelian variety ${\rm Pic}^0(S)$
is isomorphic to the ${\GG}_m$-extension of ${\rm Pic}^0(C)$
with extension class $D_1-D_2$. This ${\GG}_m$-extension $G$ 
can be realized by considering the line bundle  $L$  on ${\rm Pic}^0(C)$ 
associated to the divisor class of $D_1-D_2$ and deleting 
the zero section. It is an algebraic group since it can be 
identified with the theta group of $L$, cf.\ \cite{M,vdGM}. 

Just as in the preceding section there are two ways for compactifying
$G$: one by considering the moduli of rank $1$ torsion-free sheaves
on $S$ and secondly by glueing two sections of the ${\PP}^1$-bundle
defined by $G$ (the rank-1-compactification). The result is the same. 

We consider the corresponding ${\PP}^1$-bundle
$q: P={\PP}(L\oplus O) \to{\rm Pic}^0(C)$. 
The ${\PP}^1$-bundle $P$  contains two
effective divisors $T_1$ and $T_2$ given by ${\PP}(L \oplus 0 )$ and 
${\PP}(0\oplus O)$. There exist two sections $t_i: {\rm Pic}^0(C)
\to P$ ($i=1,2$) of $q$ with image $T_i$.  Then $t_1$ corresponds
to the projection $L\oplus O \to L$ and $t_2$ to $L\oplus O \to O$.
Since $O(T_i)\otimes O(1)^{-1}$ is trivial on the fibers of $q$
it is the pull back of a line bundle $\Lambda_i$ on ${\rm Pic}^0(C)$
and one determines $\Lambda_1\cong O$ by pulling $O(T_1)\otimes O(1)^{-1}$
under $t_2$, and similarly $\Lambda_2\cong L^{-1}$. We thus get as in 
equation (\ref{T1T2})
\begin{equation}\label{T1T2two}
O(T_1)=O(1), \qquad O(T_2)\cong O(1)\otimes q^*(L)^{-1}.
\end{equation}
We construct a non-normal variety $G^c$
by glueing $T_1$ with $T_2$ by a translation over $D_1-D_2$ 
in ${\rm Pic}^0(C)$.  It contains $G$ as a open
subvariety and the singular locus $\Sigma$ is isomorphic to 
${\rm Pic}^0(C)$.

We may interpret $G$ alternatively as a ${\GG}_m$-extension of
${\rm Pic}^0(\spr C)={\rm Pic}^0(S)$. 
Then we have the following interpretation for the compactification 
$G^c$ obtained here, cf.\ \cite{AK2}, Section 3. 
We consider the moduli space of rank $1$ torsion-free sheaves 
of $O_{S}$-modules on $S$  of first Chern class $0$. 
If $N$ is such a sheaf then outside the 
singular locus ${\rm sing}(S)$ the sheaf $N$ is locally free. 
Along ${\rm sing}(S)$ we have that $N\cong O_{S,{\rm sing}(S)}$ 
or $N\cong \nu_*(O_{\tilde{S}})|{\rm sing}(S)$.
 
Let $c=3x-\delta/2$ be the class (for algebraic equivalence) of
$C_1$ and $C_2$ on $\tilde{S}=\spr C$. Recall that $\nu$
is the normalization map $\nu: \tilde{S} \to S$. For a line bundle
$N$ in ${\rm Pic}^{-c}(\spr C)$ the direct image
$\nu_*(N)$ is a torsion free sheaf of rank $1$ of first Chern class
$0$ which is not locally free. So our situation is very similar
to the one of one-nodal curves and we have a morphism
$$
\nu_* : {\rm Pic}^{-c}(\spr C)\to \Sigma 
$$
that is an isomorphism. 
We have a natural action of $G$ on $G^c$ extending the action on itself.
On $\Sigma$ this action is compatible with the action of
${\rm Pic}^0(\tilde{S})\cong {\rm Pic}^0(C)$ on 
${\rm Pic}^{-c}(\spr C)$

We also have Abel-Jacobi maps here. If we pick a base point $p_0+q_0$
on $\spr C$ we have the map
$$
u:\spr C  \to {\rm Pic}^0(C), \qquad p+q \mapsto O(p_0+q_0-p-q).
$$
As in the case of the compactified Jacobian there is a lift of $u$ to a map $\tilde{u}:
\spr C\to P$ given by a surjection $u^*(L\oplus O)
\to M$ with $M$ an invertible sheaf on $\spr C$. In fact,
take $M=O(C_1)=O(C_2)\otimes L$. Then we have 
$$
M \otimes u^*(L\oplus O)^{\vee} =O(C_1)\oplus O(C_2)
$$
and this has a canonical section $1\oplus 1$ defining 
$u^*(L\oplus O) \to M$. This choice of $M$ is
dictated by the fact that we want $\tilde{u}^{-1}(T_i)=C_i$,
and $u^*(L)=O(C_1-C_2)$ and moreover that $\tilde{u}^*(O(1))$
should be equal to $M$. Note that the restriction of $\tilde{u}$ 
on $\spr C -C_1-C_2$ is the map defined by the canonical rational 
section $1$ of the divisor $C_1-C_2$. Moreover,
if $p_1+q_1$, $p_2+q_2$ are complementary points on the curves $C_1$, $C_2$ respectively, see Section 
\ref{divisors on symC}, then $q\tilde{u}(p_1+q_1)-q\tilde{u}(p_2+q_2) =u(p_1+q_1-p_2-q_2)=O(D_1-D_2)$. 
This implies that the map 
$\tilde{u}$ descends to an Abel-Jacobi map $\bar{u}: S \to G^c$.

Finally we calculate the class of $\tilde{u}(\spr C)$
in $P$ modulo algebraic equivalence.

\begin{proposition}
\label{classsymm2}
The algebraic equivalence class of $\gamma=\tilde{u}(\spr C)$
in $P$ is given by
$$
\gamma = q^*[C]+\frac{1}{2} q^*[C \ast C] \cdot \eta,
$$
where $C\ast C$ is the Pontryagin product and $\eta$ the class of
$O(1)$.
\end{proposition}
\begin{proof}
As above we have $[T_1]=\eta$ and since $L=O(D_1-D_2)$ it is
algebraically equivalent to $0$, hence by relation (\ref{T1T2two})
we have $[T_2]=\eta$ too. We can write our class
$\gamma$ as $q^*(a_1)+q^*(a_2)\cdot \eta$
with $a_i$ a dimension $i$ cycle on ${\rm Pic}^0(C)$.
We have $a_2=q_*(\gamma)=q_*\tilde{u}_* 1_{\spr C}=u_* 1_{\spr C}=
\frac{1}{2} [C \ast C]$ and $q_*(\gamma \cdot \eta)=[C]$.
The result follows.
\end{proof}

\end{section}
\begin{section}{The Clemens-Griffiths Map}
For smooth Fano surfaces it is natural to consider the map
$$
S \to {\rm Pic}^0(S), \qquad s \mapsto O(D_s-D_{s_0}),
$$
where $s_0$ is a fixed base point in $S$, see \cite{CG}.  
This map embeds $S$ into  ${\rm Pic}^0(S)$. There is an analogue 
of this map for the singular Fano surface $S$. We choose 
$p_0+q_0\in\spr C -C_1-C_2$ and let 
$s_0 =\nu(p_0+q_0) \in S-{\rm sing}(S)$. 
We consider the incidence variety 
$I=\{ (s,t) \mbox { with } t\in D_s\} \subset S\times S$. 
If $\pi_i: S \times S \to S$ is the $i$-projection then  
$I-\pi_1^*D_{s_0}\subset S\times S \stackrel{\pi_2}{\to} S$ 
is a family of Cartier divisors over $S-{\rm sing}(S)$, 
see Proposition \ref{nu*Dpq},  
and therefore this defines a map 
$u_0: S-{\rm sing}(S)\to {\rm Pic}^0(S)$.  

\begin{lemma}
\label{le:ajcg} Let $q: {\rm Pic}^0(S)=G \to {\rm Pic}^0(C)$ 
be the natural projection and $\nu:\spr C  \to S  $ the normalization
map. Then we have the equality of maps
$$
q \circ u_0 \circ \nu  = u :\spr C  -C_1-C_2 \to {\rm Pic}^0(C).
$$
\end{lemma}
\begin{proof}
Let $i: {\rm Pic}^0(\spr C) \to {\rm Pic}^0(C)$ 
be the natural isomorphism given in \ref{Picsym}. Then by
Corollary \ref{pullback lb} we have
$$
q \circ u_0 \circ \nu (p+q) =i\nu^*(D_{\nu(p+q)}-D_{s_0})= 
i(D_{p+q}-D_{p_0+q_0}).
$$
By relation (\ref{O(D(p+q))}) we have
$$
D_{p+q}-D_{p_0+q_0} \sim [S_{K-p-q}-\Delta/2]-[S_{K-p_0-q_0}-\Delta/2]
\sim S_{p_0+q_0-p-q}.
$$
and the result follows by the definition of $i$. 
\end{proof}

The above lemma basically says that $u_0$ is a lift to 
$G={\rm Pic}^0(S)$ of the usual 
Abel-Jacobi map $u:\spr C  \to {\rm Pic}^0(C)$. 
The next proposition shows that it coincides with the 
generalized  Abel-Jacobi map $\bar{u}$.  

\begin{proposition}
\label{ajcg}
Let  
$u_0: S-{\rm sing}(S) \to {\rm Pic}^0(S)=G$ with $u_0(s)=O(D_s-D_{s_0})$ 
be the Clemens-Griffiths map for the singular Fano surface. 
Then $u_0$ coincides with $\bar{u}$ on $S-{\rm sing}(S)$, with 
$\bar{u}$ the Abel-Jacobi map defined in section \ref{compfano}. 
\end{proposition}
\begin{proof}
Recall that the restriction to $\spr C -C_1-C_2$
of the lifting $\tilde{u}$ of the Abel-Jacobi
map $u$ is given by the canonical rational section of $O(C_1-C_2)$
on $\spr C$, as we saw in the preceding section. 
For $p+q \in\spr C $ not on $C_1$ nor on $C_2$ and $s=\nu(p+q)$ 
we have seen that the divisor 
$D_{p+q}-D_{p_0+q_0}$ on $\spr C$  descends 
to the Cartier divisor $D_s-D_{s_0}$ on $S$. 
The image $u_0(s)$ for $s\in S-{\rm sing}(S)$
is given by the class of the pull back $L=\nu^*O(D_s-D_{s_0})$
and an isomorphism $\gamma_1^*L \cong \gamma_2^*L$ on $C$.
Let $1$ be the canonical rational section of 
 $O(D_{p+q}-D_{p_0+q_0})$.  The pull back $L$ is isomorphic to $O(D_{p+q}-D_{p_0+q_0})$ and the
glueing is given by the ratio $\gamma_2^*(1)/\gamma_1^*(1)$ of the two sections. 
To determine now the embedding $u_0 : S-{\rm sing}(S) \to G={\rm Pic}^0(S)$ 
we have to carry out the above construction for the family of Cartier divisors 
$I-\pi_1^*D_{s_0}\subset S\times S \stackrel{\pi_2}{\to} S$ over $S-{\rm sing}(S)$ which defines the map $u_0$. 

Let  
${\mathcal D} = \{ (p+q,r+s) \mbox{ with }  p+q\in D_{r+s}\} 
\subset\spr C \times\spr C $.
Note that $(\nu \times \nu)_*{\mathcal D}=I$, see Proposition  \ref{nu*Dpq}. Take 
$1 \times \gamma_i :\spr C  \times C \to\spr C \times\spr C $ and we have 
$$ 
(1 \times \gamma_i)^*({\mathcal D}-\pi_1^*D_{p_0+q_0})|\{p+q\}\times C = 
p+q-p_0-q_0,  \mbox{ see Inters } \ref{Dpq dot Ci} 
$$
and 
$$  
(1 \times \gamma_i)^*({\mathcal D}-\pi_1^*D_{p_0+q_0})|\spr C \times \{p \} = 
 D_{p'_i+p''_i}- D_{p_0+q_0} \mbox{ with } p'_i+p''_i=\gamma_i(p).
$$
Let $1$ be the canonical rational section of  $O({\mathcal D})$. Then  
$\gamma^*_2(1) /\gamma^*_1(1)$ is up to a non-zero multiplicative scalar
the canonical rational section $1$ of $O(A)$ with $A$ the divisor
$$
A= (1 \times \gamma_2)^*({\mathcal D}-\pi_1^*D_{p_0+q_0}) -
(1 \times \gamma_1)^*({\mathcal D}-\pi_1^*D_{p_0+q_0})
$$ 
on $\spr C \times C$. But $A|\{p+q\}\times C $ is the zero divisor 
for every $p+q \in\spr C $. 
Hence $A$ is the pull back from $\spr C$ of the divisor 
$A|\spr C \times \{p \}=D_{p'_2+p''_2}-D_{p'_1+p''_1} = C_1-C_2$, 
see Lemma \ref{Dpq, p+q in Ci}. 
Therefore the section $\gamma^*_2(1) /\gamma^*_1(1)$ which gives 
the glueing over $S-{\rm sing}(S)$ is (up to a non-zero scalar) the pull back 
of the canonical rational section of $O(C_1-C_2)$ and hence the result. 
\end{proof}

\end{section}
\begin{section}{The Limit of the Clemens-Griffiths Map}
\label{familyfano}
Assume that we have a family $ {\mathcal X} \to \Delta $, with 
$\Delta$ the spectrum of a discrete valuation ring 
(or an open unit disc in the complex case) with generic fibre
$X_{\eta}$, a  smooth cubic threefold and special fibre $X_0$,
 a nodal cubic threefold. 
Let ${\mathcal S} \to \Delta$ be the corresponding family 
of Fano  surfaces with ${\mathcal S}_0$ the non-normal Fano surface of $X_0$. 
We may assume that the family ${\mathcal S} \to \Delta$ has a section 
$\sigma: \Delta \to {\mathcal S}$ with $\sigma (0) \in S_0-{\rm sing}(S_0)$. 

The map $S_{\eta}\times S_{\eta} \to {\rm Pic}^0 S_{\eta}$ given by
$(s,s') \to [D_s-D_{s'}]$  has generic degree $6$ and has as image a
divisor $\Theta$ that defines a principal polarization.

We consider the correspondence 
$I_{\eta}$ on $S_{\eta}\times_{{\Delta}_{\eta}} S_{\eta}$
given by pairs $(s_1,s_2)$ with $s_2 \in D_{s_1}$. This gives us
a relatively effective divisor ${\mathcal D}$ over $S_{\eta}$ via
the projection on the second factor $S_{\eta}$. In turn this defines
an embedding $\phi_{\eta}: S_{\eta} \to {\rm Pic}^0(S_{\eta}/\Delta_{\eta})$
that sends $s$ to $D_s-D_{\sigma}$.
The family $\phi_{\eta}$ is a flat family as it is irreducible and the
base $\Delta_{\eta}$ is $1$-dimensional. We consider the 
rank-$1$-compactification of the relative Picard variety with special
fibre $G^c$ and let $F$ be the flat limit of $\phi_{\eta}(S_{\eta})$ 
in the special fibre.

\begin{proposition}
The flat limit $F$ of the Fano surface coincides with the Abel-Jacobi image
$\bar{u}(S_0)$.
\end{proposition}
\begin{proof}
In characteristic $0$ it is well-known that the cohomology class of the fibre 
$\phi_{\eta}(S_{\eta})$ is $\theta^3/3!$, with $\theta$ the 
polarization class on the Picard variety. In positive characteristic
we can lift the hypersurface 
$X$ to characteristic $0$ and deduce the result from this.
Hence the cohomology class of the limit 
$F\subset G^c$ is $\theta _0 ^3/3!$ where $\theta_0$ is the limit 
of the polarization. If $\tau: P \to G^c$ is the normalization map 
then the class of $\tau ^* \theta _0$ is equal to $q^*\xi+\eta$, 
where $\xi $ is the polarization on ${\rm Pic}^0(C)$ and $\eta $
as in Proposition \ref{classsymm2}, cf.\ \cite{M}. 
Note that $(\tau^*\theta_0)^3/3!=q^*\xi^3/3!+q^*\xi^2/2!\cdot \eta$
because $\eta^2=0$ in cohomology,  which is exactly the 
cohomology class of $\bar{u}(S_0)$, see Proposition \ref{classsymm2}.

Let ${\mathcal S}^*$ be the subscheme of ${\mathcal S}$ 
obtained by removing the singular
points of $S_0$. The correspondence $I_{\eta}$ extends to $I^*$
on ${\mathcal S}^*$ in a natural way by adding the points
$(s_1,s_2) \in S_0\times S_0$ with $s_1 \in D_{s_2}$ and the map
$\phi_{\eta}$ extends also naturally to $\phi^*:{\mathcal  S}^* \to 
{\rm Pic}^0({\mathcal S}/{\Delta})$ using the section $\sigma$.
Then $\phi_{\eta}(S_{\eta}) \subset \phi^*({\mathcal S}^*)$ and hence $\phi^*(S_0^*)$
is contained in the limit $F$. Note that $S_0^*=S_0-{\rm sing}(S_0)$ and
$\phi^*(S_0^*)=u_0(S_0-{\rm sing}(S_0)=\bar{u}(S_0-{\rm sing}(S_0))$,
see Proposition \ref{ajcg}.
Hence $\bar{u}(S_0) $ is contained in $F$ and hence is a component of $F$.
But since $F$ and $\bar{u}(S_0) $ are effective cycles and have the same homology class they should be equal, as the intersection number of $\theta_0^3$ with
$F-\bar{u}(S_0)$ otherwise would be positive since $\theta_0$ is ample.
See also \cite{De} sections 7 and 8.

\end{proof}
\noindent
{\bf Acknowledgement} The second author thanks the Korteweg-de Vries Instituut
van de Universiteit van Amsterdam,  where part of this work was done, for its
support and hospitality.
\end{section}

\end{document}